\documentclass[12pt,a4paper]{amsart} 
\usepackage{amsfonts}
\usepackage{amsthm}
\usepackage{amsmath}
\usepackage{amscd}
\usepackage[latin2]{inputenc}
\usepackage{t1enc}
\usepackage[mathscr]{eucal}
\usepackage{indentfirst}
\usepackage{graphicx}
\usepackage{graphics}
\usepackage{pict2e}
\usepackage{epic}
\numberwithin{equation}{section}
\usepackage[margin=2.9cm]{geometry}
\usepackage{color,url}

\theoremstyle{plain}
\newtheorem{thm}{Theorem}[section]

\newtheorem{prop}[thm]{Proposition}

\newtheorem{conjecture}[thm]{Conjecture}

\newtheorem{problem}[thm]{Problem}
\newtheorem{example}[thm]{Example}

\begin{document}

\title{Chromatic roots and limits of dense graphs}
\thanks{The first two authors   are partially supported by ERC Consolidator Grant 648017. The first, second and fourth author is partially supported by  MTA R\'enyi
"Lend\"ulet" Groups and Graphs Research Group, and  by  the Hungarian
National Research, Development and Innovation Office -- NKFIH, OTKA
grant no.\ K109684, and also  K104206 for the second author. }

\author[P. Csikv\'ari]{P\'eter Csikv\'ari}
\address{Massachusetts Institute of Technology \\ Department of Mathematics \\
Cambridge MA 02139 \&  E\"{o}tv\"{o}s Lor\'{a}nd University \\ Department of Computer 
Science \\ H-1117 Budapest
\\ P\'{a}zm\'{a}ny P\'{e}ter s\'{e}t\'{a}ny 1/C \\ Hungary} 

\email{peter.csikvari@gmail.com}

\thanks{\emph{PC:} Partially supported by the National Science Foundation under grant no.\ DMS-1500219 and  by the  ELTE-MTA Geometric and Algebraic Combinatorics Research Group.}

\author[P.E. Frenkel]{P\'eter E.\ Frenkel}
\address{Alfr\'ed R\'enyi Institute of Mathematics \\
H-1053 Budapest \\
Re\'altanoda u. 13-15. \\
Hungary \& E\"{o}tv\"{o}s Lor\'{a}nd University \\ Department of Algebra and
  Number Theory \\ H-1117 Budapest
\\ P\'{a}zm\'{a}ny P\'{e}ter s\'{e}t\'{a}ny 1/C \\ Hungary} 

\email{frenkelp@cs.elte.hu}

\author[J. Hladk\'y]{Jan Hladk\'y}
\address{The Czech Academy of Sciences, Institute of Computer Science, Pod Vod\'{a}renskou v\v{e}\v{z}\'{\i} 2, 182 07 Prague, Czech Republic. With institutional support RVO:67985807.
\emph{Address at time of submission:} Institute of Mathematics, Czech Academy of Sciences. \v Zitn\'a 25, 110 00, Praha, Czech Republic. The Institute of Mathematics of the Czech Academy of Sciences is supported by RVO:67985840.}\thanks{\emph{JH:} The research leading to these results has received funding from the People Programme (Marie Curie Actions) of the European Union's Seventh Framework Programme (FP7/2007-2013) under REA grant agreement number 628974. Most of the work was done while the author was an EPSRC Fellow at the Mathematics Institute, University of Warwick, UK}
\email{honzahladky@gmail.com}
\author[T. Hubai]{Tam\'as Hubai}
\address{Department of Computer Science, University of Warwick, Coventry CV4 7AL,
United Kingdom.
}
\email{t.hubai@warwick.ac.uk}
\thanks{\emph{TH:} Partially supported by EPSRC grant no.\ EP/M025365/1.}

\begin{abstract}
In this short note we observe that recent results of Ab\'ert and Hubai
and of Csikv\'ari and Frenkel about Benjamini--Schramm continuity of
the holomorphic moments of the roots of the chromatic polynomial extend
to the theory of dense graph sequences. We offer a number of problems
and conjectures motivated by this observation.
\end{abstract}

\maketitle

\section{Introduction}

Recently, there has been much work on developing limit theories of
discrete structures, and of graphs in particular. The best understood
limit concepts are those for dense graph sequences and bounded-degree graph sequences. 
The former one was developed by Borgs, Chayes, Lov\'asz, S\'os, Szegedy and Vesztergombi~\cite{BCLSV,LovSze}, and the 
latter was initiated by Benjamini and Schramm~\cite{BS}.
The convergence notion in both these theories is based on frequencies
of finite subgraphs, and it is a fundamental programme to understand
what other parameters are captured in the limits (i.e.,  are continuous
with respect to the corresponding topologies). In this short note we
show that recent proofs of Ab\'ert and Hubai and of Csikv\'ari and Frenkel
about the convergence of holomorphic moments of the chromatic roots
in a Benjamini--Schramm convergent sequence translate to the dense
model as well. Furthermore, we conjecture that in the dense
model we actually  have weak convergence of the root distributions. 

Let us now give the details. We assume the reader's familiarity with
basics of graph limits, we shall however give pointers to  Lov\'asz's recent monograph~\cite{Lov:Book} throughout. 

Recall that given a graph $G$ of order $n$, its \emph{chromatic polynomial} $P(G,x)$ (in a complex variable $x$) is defined as
\begin{equation}\label{eq:defchromatic}
P(G,x)=\sum_{k=0}^n \mathrm{ip}(G,k)x(x-1)\ldots(x-k+1)\;,
\end{equation}
where $\mathrm{ip}(G,k)$ is the number of partitions of $V(G)$ into $k$ non-empty independent sets. In other words, for nonnegative 
integer values of $x$, $P(G,x)$ counts the number of proper vertex-colorings of $G$ with $x$ colors.

Next, we recall the result of Ab\'ert and Hubai \cite{AbeHub}, which concerns convergence of
bounded-degree graphs. To this end, recall that a sequence of graphs $(G_n)_n$ of maximum
degree uniformly upper-bounded by $D$ is \it Benjamini--Schramm convergent \rm if
for each fixed connected graph $F$, the number sequence $\hom (F, G_n)/
v(G_n)$  converges. Here $\hom (F,G)$ denotes the number of homomorphisms from $F$ to $G$, and $v(G)$ denotes the number of nodes in $G$.  There are many equivalent definitions of Benjamini--Schramm convergence. A detailed treatment appears in~\cite[
Chapter 19]{Lov:Book}. 

Suppose that $G$ is a graph of maximum degree at most $D$. 
We can associate to it the uniform
probability measure $\mu_{G}$ on the multiset of the roots of
the chromatic polynomial $P(G,x)$. The Sokal bound \cite{Sokal:Bound}
tells us that this \it chromatic measure  \rm $\mu_{G}$ is supported in the disk of radius (strictly
less than) $8D$. The main result of~\cite{AbeHub} then reads as
follows.

\begin{thm} \label{thm:AbHu}
Suppose that $(G_{n})_{n}$ is a Benjamini--Schramm
convergent sequence of graphs of maximum degree at most $D$. Suppose
that $f:B\rightarrow\mathbb{C}$ is a holomorphic function defined
on the open disk $B=B(0,8D)$. Then the sequence $$\int f(z)\mathrm{d}\mu_{G_{n}}(z)$$
converges.
\end{thm}

Note that to prove Theorem~\ref{thm:AbHu} it suffices to prove the convergence
of the 
\emph{holomorphic moments} $$\int z^{k}{\mathrm{d}}\mu_{G_{n}}(z)\qquad  (k\in\mathbb{N})$$ for a  Benjamini--Schramm convergent graph sequence $(G_n)_n$, and this is indeed how the proof goes. 

As was noted in~\cite{AbeHub}, it is not always the case that the
measures $\mu_{G_{n}}$ in a Benjamini--Schramm convergent graph sequence
converge weakly. This can be seen from the following example.

\begin{example}
\label{ex:Counterexample}
Consider  paths $P_{n}$ and
cycles $C_{n}$ of growing order. These two sequences have the same
Benjamini--Schramm limit but the weak limit of $(\mu_{P_{n}})_{n\rightarrow\infty}$
is concentrated on 1 whereas the weak limit of $(\mu_{C_{n}})_{n\rightarrow\infty}$
is the uniform measure on the unit circle with the center in~1.
\end{example}

Csikv\'ari and Frenkel~\cite{CsiFre} generalized Theorem~\ref{thm:AbHu} to
a wider class of graph polynomials. This is discussed in Section~\ref{ssec:otherpoly}.\medskip{}

Let us now turn to dense graphs. The convergence notion in the dense model 
was introduced by Borgs, Chayes, Lov\'asz, S\'os, Szegedy and Vesztergombi. Of the
many equivalent definitions, we shall give the one that is the most convenient for our
purposes. We refer to \cite[Chapters 11-12]{Lov:Book} for more details. A sequence of graphs
$(G_n)_n$ is \it convergent in the dense model \rm if for each fixed graph $F$, the number sequence  $\hom (F,G_n)/
v(G_n)^{v(F)}$ converges. Let us also recall that
if we have a convergent sequence of dense graphs then we can associate to it a limit
object, a so-called \it graphon, \rm  see \cite[\S 11.3]{Lov:Book}.

 Suppose that $G$ is a graph of order
$n$. Then the vertices of $G$ have arbitrary degrees between
$0$ and $n-1$. The measure $\mu_{G}$ need not be supported in a
bounded region for a such a graph $G$; the Sokal bound gives
only a bound of roughly $8n$ on the modulus of the chromatic roots.
This bound can probably be improved down to $n-1$ (see Conjecture~\ref{conj:modulusN}) but not more.
Thus, it is natural to scale down $\mu_{G}$ by the factor of $n$, defining
a new probability measure $\nu_{G}$, $\nu_{G}(X):=\mu_{G}(nX)$, where for $X\subset\mathbb{C}$ we define $nX=\{nx:x\in X\}\subset\mathbb{C}$.
Now, $\nu_{G}$ is supported in the disk of radius $8$. The main
result of this note is the observation that Theorem~\ref{thm:AbHu} has a counterpart
for sequence of dense graphs.

\begin{thm} \label{thm:Dense} Suppose that $(G_{n})_{n}$ is a sequence of graphs which converges in the dense model. Suppose that $f:B\rightarrow\mathbb{C}$
is a holomorphic function defined on an open disk $B=B(0,8)$. Then
the sequence $$\int f(z)\mathrm{d}\nu_{G_{n}}(z)$$ converges.
\end{thm}

We will give a sketch of a proof of Theorem~\ref{thm:Dense} in Section~\ref{sec:ProofThmDense}.

\medskip

Note that by a standard argument from complex analysis, we can approximately count the number of colorings in a convergent graph sequence. 
Note that when $G$ has $n$ vertices and $\ell=Cn$, we expect $P(G,\ell)$ to grow as $(cn)^n$ for some $c\in\mathbb R$.

\begin{thm}
Let $(G_n)_n$ be a sequence of graphs convergent in the dense model, where $G_n$ has order $n$.  Then for each $C>8$, the quantity
$$\frac{\sqrt[n]{P(G_n,Cn)}}n$$ converges as $n\rightarrow \infty$.
\end{thm}

The proof follows the same lines as Theorem 1.2 of ~\cite{AbeHub}. Also, the result can be stated a bit more generally, as can again be seen in Theorem 1.2 of ~\cite{AbeHub}.

\medskip
We believe that there exists no counterpart to Example~\ref{ex:Counterexample}
for dense graphs. This is the main conjecture of the present paper.

\begin{conjecture}
\label{conj:weak}Suppose that $(G_{n})_{n}$ is a sequence of graphs convergent in the dense model. Then the rescaled chromatic measures  $\nu_{G_{n}}$ converge
weakly.
\end{conjecture}

In general, a graphon does not carry much information about chromatic
properties of graphs which converge to it. For example, it is easy
to construct a sequence of graphs such that their chromatic numbers grow
almost linearly with their orders, yet converge to the constant-zero
graphon.  On the other hand it is easy to construct another sequence of graphs such that their chromatic numbers grow arbitrarily slowly, yet converge to the constant-one graphon.  That is, in a sense, the chromatic number is not even semicontinuous with respect to the cut-distance.

An immediate consequence of Conjecture~\ref{conj:weak}
would be that it would allow us to associate ``chromatic roots''
to graphons. This is perhaps the most substantial information about
chromatic properties which could be reflected in the limit.

The only support for Conjecture~\ref{conj:weak} is a lack of counterexamples
we could come up with. In particular, the Conjecture asserts that
the normalized chromatic measures of Erd\H{o}s--R\'enyi random graphs
(with constant edge probability) or more generally random graphs coming
from sampling from a graphon converge --- and this seems to be a very
weak form of the conjecture. It would be very interesting to prove
this, and to describe the weak limit.

\begin{problem}
\label{prob:ErdosRenyi}
What is the typical distribution of the chromatic
roots of the Erd\H{o}s--R\'enyi random graph $\mathbb{G}_{n,p}$, for a fixed
$p\in(0,1)$?
\end{problem}

Computational restrictions allowed us to run simulations only for $n\le 10$. Such limited simulations did not hint for any limit behavior.

\medskip
Last, let us remark that the measure $\nu_{G}$ is not trivialized
by the scaling we introduced. This is stated in the following proposition.

\begin{prop} \label{prop:nottrivial}
For every $\delta>0$ there exists $\epsilon>0$
such that the following holds. Suppose that $G$ is a graph of order
$n$ with at least $\delta n^{2}$ edges. Then at least $\epsilon n$
of the chromatic roots of $G$ have modulus at least $\epsilon n$.
\end{prop}

\section{\label{sec:ProofThmDense}Proof of Theorem~\ref{thm:Dense}}

It is only needed to observe that the argument in~\cite{AbeHub}
is valid even in the dense model. More precisely, in \cite[Theorem 3.4]{AbeHub}
the following is proven. The symbol $\hom(T,H)$ stands for the number
of homomorphisms from a graph $T$ to a graph~$H$.
\begin{thm}[{\cite[Theorem 3.4]{AbeHub}}]
Let $H$ be a graph, and for $k\in\mathbb{N}$ let $$p_{k}=|V(H)|\int z^{k}\mathrm{d}\mu_{H}(z).$$ Then
\begin{equation}
p_{k}=\sum_{T}(-1)^{k-1}kc_{k}(T)\hom(T,H)\;,\label{eq:pk}
\end{equation}
where $c_{k}(T)$ are constants, and the summation ranges over connected
graphs $T$ of order at most $k+1$.
\end{thm}
With this result, the proof of Theorem~\ref{thm:Dense} is straightforward. Let us write $p_{k,n}$ for the number $p_{k}$ from the previous
theorem associated to the graph $G_{n}$. As was remarked earlier,
it suffices to prove the theorem for $f(z)=z^{k}$, $k\in\mathbb{N}$.
For simplicity, let us assume that the graph $G_{n}$ has $n$ vertices.
We have
\[
\int z^{k}\mathrm{d}\nu_{G_{n}}(z)=\frac{1}{n^{k}}\int z^{k}\mathrm{d}\mu_{G_{n}}(z)=\frac{p_{k,n}}{n^{k+1}}\;.
\]
The sequence $(G_{n})$ is convergent. In particular, for every graph
$T$ of order at most $k+1$ the quantity $$\frac{\hom(T,G_{n})}{n^{k+1}}$$
converges. Observe that the right-hand side of (\ref{eq:pk}) (for
a fixed number $k$) contains only a bounded number of summands. Consequently,
$$\frac{p_{k,n}}{n^{k+1}}$$ converges as $n\rightarrow\infty$, finishing
the proof.

\section{Proof of Proposition~\ref{prop:nottrivial}}

%We use~(\ref{eq:pk}) again. The table at the end of~\cite{AbeHub} gives us the constant $c_{1}(\mbox{edge})$ from the right-hand side of~\eqref{eq:pk},  $c_{1}(\mbox{edge})=\frac{1}{2}$. That is, 
%\[
%\frac{1}{2}\hom(\mbox{edge\ensuremath{,H)}}=p_{1}(H)=n\int z\mbox{d}\mu_{G_{n}}(z)=\sum_{x\;\mathrm{ chr.root}}x\;.
%\]

It is well known, and easy to see from the formula \eqref{eq:defchromatic}, that
the sum of the chromatic roots of $G$ is the number of edges in $G$.
By the assumption of the proposition, this
%the left-hand side
 is at least
$\delta n^{2}$. Also, recall that the chromatic roots are contained in the disk of radius~$8n$. Thus, for $$\delta n^2\le \sum_{x\;\mathrm{ chr.root}}x$$
to hold, we must have  at least $\frac{\delta}{9}n$ roots $x$
of the chromatic polynomial with $\Re(x)\in\left[\frac{\delta}{9}n,8n\right]$.

\section{Remarks and conjectures}

\subsection{Variants of Sokal's bound}
Recall that the bound asserts
that if a graph has maximum degree $\Delta$, then all the chromatic
roots lie in the disk of radius $r=8\Delta$. The value $8\Delta$
is not optimal; Sokal himself actually gives $7.96\ldots\times \Delta$.
On the other hand the complete bipartite graph $K_{\Delta,\Delta}$
shows~\cite{Sokal:Dense}
that one cannot go below $r=1.59\ldots\Delta$.\footnote{Let us note that this has not been proven rigorously.} Here, we suggest to bound the moduli of the chromatic roots by the
order instead of the maximum degree.

\begin{conjecture} \label{conj:modulusN}
Every graph $G$ of order $n$ has all the chromatic
roots of modulus at most $n-1$.
\end{conjecture}

If true, complete graphs would be the extremal graphs for the problem.\footnote{Recall that the roots of the chromatic polynomial of a complete graph $K_n$ are $\{0,1,\ldots,n-1\}$.} Note that  Conjecture~\ref{conj:modulusN} is known to be true for real zeros. Indeed, if $x>n-1$ is real then each summand in~\eqref{eq:defchromatic} is non-negative, and the summand for $k=n$ is strictly positive, yielding $P(G,x)>0$. Secondly, we claim that if $x$ is negative then it is not a root of $P(G,\cdot)$. Indeed, it is well-known (see e.g.~\cite[Corollary~2.3.1]{Dong:Chromaticbook}) that the coefficients of $P(G,\cdot)$ alternate in sign. The value $P(G,x)$ is then a sum of terms with the same sign, and in particular, non-zero.

By enumerating all graphs of a given order on a computer, we have verified Conjecture~\ref{conj:modulusN} for $n\le 10$.

\medskip
Our next problem can be seen as an extension of Sokal's bound, but is also connected to Conjecture~\ref{conj:weak} as we show below.
\begin{problem} \label{prob:addedges}
Suppose that $G$ is a graph and $G'$ is obtained
from $G$ by adding edges in such a way that the degree at each vertex
increases by at most $\Delta$. Is it true that the chromatic roots
move by at most $c\Delta$, for some absolute constant $c$? (By ``moving''
we mean that there is a bijection $\pi$ from the multiset of the
chromatic roots of $G$ to the multiset of the chromatic roots of
$G'$ so that $\left|x-\pi(x)\right|\le c\Delta$ for each chromatic
root $x$ of $G$.)
\end{problem}
Note that to answer Problem~\ref{prob:addedges} in the affirmative, it
would suffice to prove the case $\Delta=1$.

Suppose that $G_1$ and $G_2$ are two $n$-vertex graphs with edit-distance at most $\epsilon n^2$. That means that after a suitable vertex identification of $V(G_1)$ and $V(G_2)$ the graph $G$ on the same vertices whose edges are the common edges of $G_1$ and $G_2$ has the property that for at most $2\sqrt{\epsilon}n$ vertices do we have $\deg_G(v)\le \deg_{G_1}(v)-\sqrt{\epsilon}n$ or  $\deg_G(v)\le \deg_{G_2}(v)-\sqrt{\epsilon}n$. For the sake of drawing the link to Conjecture~\ref{conj:weak}, let us assume that there are no such exceptional vertices. Then a positive solution to Problem~\ref{prob:addedges} would give that the chromatic measure $\nu_{G_1}$ and $\nu_{G_2}$ are close in the weak$^*$ topology. In particular, a positive answer to Problem~\ref{prob:addedges} would provide a support for Conjecture~\ref{conj:weak} when the topology generated by the cut-distance is replaced by the stronger $L^1$-metric.

\subsection{Matching polynomial}\label{ssec:otherpoly}
As mentioned before, Theorem~\ref{thm:AbHu} has been~\cite{CsiFre} extended to
a large class of ``multiplicative graph polynomials of bounded exponential
type''.  (In particular, this includes univariate polynomials derived
from the Tutte polynomial, and a modified version of the  matching
polynomial.  For the matching polynomial, the behavior of the root distribution in
Benjamini--Schramm convergent graph sequences was discussed in \cite{ACSFK,
  ACSH}.)

 The proof in~\cite{CsiFre} of this more general statement  translates
verbatim to the dense setting as well,\footnote{Let us remark that the proof
  is quite different from the original proof by Ab\'ert and Hubai.} thus
giving Theorem~\ref{thm:Dense} for multiplicative graph polynomials of bounded
exponential type. Problem~\ref{prob:ErdosRenyi} can be asked for these alternative graph polynomials as well. 

The case of the matching polynomial is particularly simple. We recall the definition.
Let $G=(V,E)$ be a finite graph on $v(G)=n$ vertices. Let $m_k(G)$ be the number of $k$-matchings. Note that $m_0(G)=1$ and $m_k(G)=0$ for $k>\lfloor n/2\rfloor$. The \emph{matching polynomial} $\mu(G,x)$ in one variable $x$ is defined as
\begin{align*}
\mu(G,x)&=\sum_{k=0}^{\lfloor n/2\rfloor}(-1)^km_k(G)x^{n-2k}.
\end{align*}
A well-known result of Heilmann and Lieb \cite{hei} asserts that the roots
of the matching polynomial are all real. 
It is easy to see that the matching
polynomial is multiplicative (w.r.t.\ disjoint union) and the coefficient of
$x^{n-i}$ is a linear combination of subgraph counts. Thus, the version of the main
result of~\cite{CsiFre} for dense graphs applies. Recall the Weierstrass approximation theorem: on a compact interval, any continuous function can be uniformly approximated by polynomials.
It follows that
convergence  of
(holomorphic) moments, i.e., $$\int x^k\mathrm d\nu_n\to\int x^k\mathrm d\nu,\quad (k=1,2,\dots)$$ is equivalent to the 
weak convergence $\nu_n\to\nu$ for probability distributions $\nu_n$, $\nu$ 
supported on a  compact interval. We get the following. If $(G_n)$ is a sequence of
graphs converging in the dense model, consider the uniform distribution
$\pi_n$ on the roots of the matching polynomial $\mu(G_n,x)$. Then %the
%rescaled distributions 
$%\frac{1}{v(G_n)}
\pi_n$ scaled down by a factor of $v(G_n)$ converges weakly. Let us explain
why this corollary is trifling. Indeed, the full Heilmann--Lieb theorem
asserts that if $G$ is a graph of maximum degree $D$, and $G$ is not a
matching, than the roots of the matching polynomial $\mu(G,x)$ lie in the
interval $\left[-2\sqrt{D-1},2\sqrt{D-1}\right]$, and in particular in
$\left[-2\sqrt{v(G)-2},2\sqrt{v(G)-2}\right]$. In other words, the %rescaled
distribution $%\frac{1}{v(G_n)}
\pi_n$ scaled down by a factor of $v(G_n)$ converges to the Dirac measure at 0. So,
the rescaling suggested by the Heilmann--Lieb theorem is by a factor of 
$%\tfrac{1}
{\sqrt{v(G_n)}}
 $. To get the right statement, we need to
introduce the \emph{modified matching polynomial}. This is a polynomial in one variable $x$ defined by
$$M(G,x)=\sum_{k=0}^{\lfloor n/2\rfloor}(-1)^km_k(G)x^{n-k}.$$ 
The matching polynomial and its modified version encode the same information. Indeed, we have $\mu(G,x)=x^{-n}M(G,x^{2})$. %Consequently,
 We can factor $M(G,x)$ as
\begin{equation}
\label{WER}
M(G,x)=x^{\lceil n/2\rceil}\prod_{i=1}^{\lfloor n/2\rfloor}(x-\gamma_i(G) )\;.
\end{equation}
Then the real numbers$$\left(\pm \sqrt{\gamma_i(G)}\right)_{i=1}^{\lfloor n/2\rfloor},$$ together with an
extra zero if $n$ is odd,  are the roots of $\mu(G,x)$. 
%(Note that this is possible even when the degree of $M(G,x)$ is less than $\lfloor n/2\rfloor$, and that this factorization is unique up to the order.)

It can be easily checked directly from~\cite[Definitions~1.3,~1.4]{CsiFre}
that $M(G,x)$ is %also
 a graph polynomial of bounded exponential type. So, it is  the \it modified
 \rm matching polynomial that we want to apply the main result
 of~\cite{CsiFre} to. We thus readily obtain a counterpart  of Conjecture~\ref{conj:weak} for the roots of the matching polynomial, with the right scaling.
 
\begin{thm}\label{thm:matchingconvergence}
Suppose that $(G_{n})_{n}$ is a
sequence of graphs convergent in the dense model. Let $\pi_n$ be the uniform probability measure on roots of the matching polynomial $\mu(G_n,x)$. Then the rescaled measures $$\lambda_n(X):=\pi_n\left({\sqrt{v(G_n)}}X\right)$$ converge
weakly.
\end{thm}
In particular, this allows us to associate a ``matching measure'' to a graphon (cf.\ text below Conjecture~\ref{conj:weak}).

\medskip

In the rest of this section we answer the counterpart of Problem~\ref{prob:ErdosRenyi} for the matching polynomial. This was done independently, and prior to the current manuscript being publicly available, in~\cite{ChLiLi:Matching}.  Our proof relates the roots of the matching polynomial of $G_n$ to those of the complete graphs $K_n$. To compare, the proof in~\cite{ChLiLi:Matching} goes via counting ``tree-like walks'', a concept introduced in~\cite{GodsilTreeWalks}.

We can now state the main result of~\cite{ChLiLi:Matching}.
\bigskip

\begin{thm}
Let $p\in(0,1)$, and let $(G_n)_n$ be a sequence of Erd\H{o}s--R\'enyi random
graphs $G_n\sim \mathbb{G}_{n,p}$. Let $\pi_n$ be the uniform probability
distribution on the roots of the matching polynomial of $G_n$. Then almost
surely, the measures $\lambda_n(X):=%\frac{1}
\pi_n\left(\sqrt n X\right)$ converge weakly to the semicircle distribution $SC_p$ whose density function is
$$\rho_p(x):=\frac1{2\pi}\sqrt{4-\tfrac{x^2}p}\;, \quad -2p\le x\le 2p\;.$$
\end{thm}
In combination with Theorem~\ref{thm:matchingconvergence}, this determines the limit of matching measures for an arbitrary sequence of quasirandom (in the sense of Chung--Graham--Wilson) graphs. Also, note that we present a proof only for $p$ fixed, but the same technique works also for $$p=\Omega\left(\frac{\log^{\text{const}}n}{n}\right).$$

\begin{proof}
Since all the roots of the matching polynomial are real, the convergence of the holomorphic moments $\int z^k\mathrm{d}\lambda_n(z)$, $k\in\mathbb{N}$ readily implies convergence in distribution. Let us thus argue that for each $k\in\mathbb{N}$, almost surely we have $$\int z^k\mathrm{d}\lambda_n(z)\rightarrow \int z^k\mathrm{d} SC_p(z)\;.$$
 For each fixed $k=0,1,2,\ldots$, and for the random graphs $G_n$, we asymptotically almost surely have
\begin{equation}\label{eq:POI}
\frac{m_k(G_n)}{m_k(K_n)}=(1+o(1))p^k\;,
\end{equation}
as each set of $k$ pairs of vertices has probability $p^k$ of being entirely
included as edges of $G_n$, and this quantity is concentrated around the
expectation. For details on  how to prove such a result, see \cite[Chapter 4]{AlonSpencer}.
Since the $m_i(G)$ are elementary symmetric polynomials of the roots $\gamma_i(G)$ (cf.~\eqref{WER}), the Newton identities give that for each fixed $k$,
\begin{equation*}
\sum_{i=1}^{\lfloor n/2\rfloor}\gamma_i(G)^k=P_k(m_1(G),\dots ,m_k(G))
\end{equation*}
for some multivariate polynomial $P_k$. It follows from the Newton identities that the polynomial $P_k$ has the property that
$$P_k(a_1t,\dots ,a_kt^k)=t^kP(a_1,\dots ,a_k).$$
Putting this together with~\eqref{eq:POI}, we get that
$$\sum_{i=1}^{\lfloor n/2\rfloor}\gamma_i(G)^k=(1+o(1))p^kP_k(m_1(K_n),\dots ,m_k(K_n)).$$
We conclude that
$$\sum_{i=1}^{\lfloor n/2\rfloor}\gamma_i(G)^k=(1+o(1))p^k\sum_{i=1}^{\lfloor n/2\rfloor}\gamma_i(K_n)^k.$$
By a classical result of Heilmann and Lieb~\cite[(3.15)]{hei}, the matching polynomials of complete graphs are the Hermite polynomials. The distribution of zeros of the Hermite polynomial of degree $n$ scaled down by $\sqrt{2n}$ converges to the semicircle distribution $SC_1$, see for instance \cite{lal}. Hence almost surely the measures $\lambda_n$ converge weakly to the semicircle distribution $SC_p$. (Note that the zeros of the matching polynomial are supported on $\pm \sqrt{\gamma_i}$ so we have to rescale the semicircle distribution only by a factor of $\sqrt{p}$.)
\end{proof}

\section{Acknowledgements}
Most of the work was done in the summer of~2013 while JH was visiting E\"otv\"os Lor\'and University. 
He would like to thank L\'aszl\'o Lov\'asz for helping him with the arrangements and all the members of the group for a stimulating atmosphere.

\bigskip

The contents of this publication reflect only the authors' views and not necessarily
the views of the European Commission or the European Union.

\end{document}